\documentclass[11pt,english]{amsart}
\usepackage[T1]{fontenc}
\usepackage[latin9]{inputenc}
\usepackage{amsbsy}
\usepackage{amstext}
\usepackage{amsthm}
\usepackage{amssymb}

\makeatletter
\usepackage{amsmath,amsfonts,amssymb,amsthm,epsfig}

\usepackage{color}
\usepackage[final,allcolors=blue,colorlinks=true]{hyperref}

\def\R{\mathbb R}

\def\S{\mathbb S}
\def\X{\mathbb X}

\def\al{\alpha}
\def\be{\beta}

\def\ep{\epsilon}

\def\na{\nabla}
\def\Om{\Omega}  
\def\De{\Delta}      

\def\wq{\infty}
\def\pa{\partial}

\def\id{\text{\rm id}}

\newcommand{\dist}{\text{\rm dist}}
\newcommand{\medint}{-\kern -,375cm\int}         
\newcommand{\medintinrigo}{-\kern -,315cm\int}

\numberwithin{equation}{section}
\newtheorem{theorem}{Theorem}[section]
\newtheorem*{theorem*}{Theorem}  

\newtheorem*{conclusion*}{Conclusin}

\newtheorem{corollary}[theorem]{Corollary}
\newtheorem*{corollary*}{Corollary}

\newtheorem*{lemma*}{Lemma}

\newtheorem*{notation*}{Notation}

\newtheorem*{proposition*}{Proposition}

\newtheorem*{remark*}{Remark}

\newtheorem*{example*}{Example}                
\theoremstyle{definition}

\makeatother

\usepackage{babel}
\begin{document}
\title[]{Conservation law of harmonic mappings in supercritical dimensions}

    \author[C.-Y. Guo and C.-L. Xiang]{Chang-Yu Guo and Chang-Lin Xiang}

\address[Chang-Yu Guo]{Research Center for Mathematics and Interdisciplinary Sciences, Shandong University, 266237,  Qingdao and  Frontiers Science Center for Nonlinear Expectations, Ministry of Education, P. R. China} \email{changyu.guo@sdu.edu.cn}
\address[Chang-Lin Xiang]{Three Gorges Mathematical Research Center, China Three Gorges University,  443002, Yichang,  P. R. China} \email{changlin.xiang@ctgu.edu.cn}
\thanks{}
\thanks{C.-Y. Guo is supported by the Young Scientist Program of the Ministry of Science and Technology
	of China (No. 2021YFA1002200), the National Natural Science Foundation of China (No. 12101362), the Taishan Scholar Project and the Natural Science Foundation of Shandong Province (No. ZR2022YQ01, ZR2021QA003). C.-L. Xiang is financially supported by the National Natural Science Foundation of China (No. 12271296).}

\begin{abstract}
In this short note, we provide a partial extension of Rivi\`ere's convervation law in higher dimensions under certain Lorentz integrability condition for the connection matrix. As an application, we obtain a conservation law for weakly harmonic mappings around regular points in supercritical dimensions.
\end{abstract}

\maketitle

{\small
\keywords {\noindent {\bf Keywords:} Conservation law, Harmonic mappings, Gauge transform, Lorentz-Sobolev spaces, Supercritical dimension}
\smallskip
\newline
\subjclass{\noindent {\bf 2020 Mathematics Subject Classification:} 58E20, 35J60}
}
\bigskip

\section{Introduction and main result}

Let $B^{n}\subset\R^{n}$ ($n\ge2$)
be the unit open ball centered at the origin and $N^{k}$ a closed $k$-dimensional $C^{2}$-Riemannian manifold which is isometrically embedded into $\R^{m}$.
A weakly harmonic mapping from $B^{n}$ into $N$ is defined
as a mapping $u\in W^{1,2}(B^{n},N)$, which is a critical point of
the Dirichlet energy $\int_{B^{n}}|\na u|^{2}$ with respect to the
outer variation. Its Euler-Lagrange equation can be written as
\begin{equation}
-\De u=A(u)(\na u,\na u),\label{eq: HM}
\end{equation}
where $A$ is the second foundamental form of $N\hookrightarrow\R^m$. A basic problem
on the theory of weakly harmonic mappings is to find the optimal regularity.

In the case when $N=\S^{m-1}\subset\R^{m}$ is the standard Euclidean sphere, several authors (\cite{Chen-1989,KRS-1989,Shatah-1988}) have independently discovered that the following system (consisting of divergence-free vector fields)
\begin{equation}
{\rm div}(u^{i}\na u^{j}-u^{j}\na u^{i})=0,\qquad\text{for all }\,1\le i,j\le m,\label{eq: Shatah-Chen}
\end{equation}
are indeed equivalent with \eqref{eq: HM}.
Since this system has a nice interpretation by Noether's conservation law,
equation \eqref{eq: Shatah-Chen} is called the conservation law of (sphere-valued)
harmonic mappings. This conservation law is very useful in investigating the regularity
and compactness theory of (sphere-valued) harmonic mappings; see e.g.~H\'elein \cite{Helein-2002} and the references therein fore more information.

In view of the above discovery and applications, Rivi\`ere asked, in his seminar work \cite[Page 4]{Riviere-2007}, the following fundamental question: \textbf{how to write the harmonic mapping equations \eqref{eq: HM} into divergence form for all dimensions $n\ge2$ and for all closed manifolds $N$?}  He succeeded to do this in \cite{Riviere-2007} when $n=2$.
The idea is to write \eqref{eq: HM} in a more general form
\begin{equation}
-\De u=\Om\cdot\na u\label{eq: Riviere}
\end{equation}
for some $\Om\in L^{2}(B^{n},so_{m}\otimes\R^{n})$. Indeed, suppose
$u$ is a smooth harmonic mapping from $B^{n}$ into $N$, and let
$\left\{\boldsymbol{\nu}_{I}=(\nu_{I}^{1},\cdots,\nu_{I}^{m})\right\}_{I=1}^{m-k}$
($I=1,\cdots,m-k$) be an orthonormal frame of $T^{\bot}N$ in a neighborhood
of $u(x)$. Then
\[
A(\na u,\na u)=-\sum_{I}\left\langle \pa_{\al}u,\na_{\pa_{\al}u}\nu_{I}\right\rangle \nu_{I}=-\sum_{I}\left\langle \pa_{\al}u,\pa_{\al}(\nu_{I}(u))\right\rangle \nu_{I}.
\]
Note that $\langle\pa_{\al}u,\nu_{I}(u)\rangle=0$ for all $1\le\al\le n$
and $I\le m-k$. Thus, in coordinates, there holds
\[
-\De u^{i}=-\sum_{j,l,I}\pa_{\al}u^{l}(\pa_{j}\nu_{I}^{l}(u)\pa_{\al}u^{j})\nu_{I}^{i}(u)=\Om_{j}^{i}(u,\na u)\cdot\na u^{j},\qquad1\le i\le m,
\]
where
\[
\Om_{j}^{i}(u,\na u)=\sum_{l,I}\left(\nu_{I}^{i}(u)\pa_{j}\nu_{I}^{l}(u)-\nu_{I}^{j}(u)\pa_{i}\nu_{I}^{l}(u)\right)\na u^{l}
\]
is antisymmetric with respect to $(i,j)$ and grows as $O(|\na u|)$.
Hence $\Om\in L^{2}(B^{n},so_{m}\otimes\R^{n})$.

Working with equation \eqref{eq: Riviere}, Rivi\`ere \cite{Riviere-2007} proposed the following theorem concerning conservation law.

\begin{theorem}[Theorem I.3, \cite{Riviere-2007}]\label{thm:Riviere}
Fix $m\in\mathbb{N}$ and $\Omega=(\Omega_{j}^{i})_{1\leq i,j\leq m}\in L^{2}\left(B^{n},so_{m}\otimes\wedge^{1}\mathbb{R}^{n}\right)$. If there are $A\in L^{\infty}\left(B^{n},M_{m}(\mathbb{R})\right)\cap W^{1,2}$
and $B\in W^{1,2}\left(B^{n},M_{m}(\mathbb{R})\otimes\wedge^{2}\mathbb{R}^{n}\right)$
satisfying
\begin{equation}\label{eq:for A and B}
d_{\Omega}A:=dA-A\Omega=-d^{*}B,
\end{equation}
then every solution to \eqref{eq: Riviere} on $B^{n}$ satisfies the following conservation law
\begin{equation}\label{eq: CL of Rivere}
d\left(*Adu+(-1)^{n-1}(*B)\wedge du\right)=0.
\end{equation}
\end{theorem}

In coordiantes, \eqref{eq:for A and B} means $dA_{j}^{i}-\sum_{k=1}^{m}A_{k}^{i}\Omega_{j}^{k}=-d^{*}B_{j}^{i}$
for all $i,j\in\{1,\ldots,m\}$ and \eqref{eq: CL of Rivere} reads as
\[
\frac{\pa}{\pa x_{\al}}\left(A_{j}^{i}\frac{\pa u^{j}}{\pa x_{\al}}+(-1)^{n-1}B_{j\al\be}^{i}\frac{\pa u^{j}}{\pa x_{\be}}\right)=0,
\]
where the Einstein summation convention is used.

Thus according to Theorem \ref{thm:Riviere}, in order to find the conservation law for harmonic mappings, it suffices to prove the existence of $A,B$ for a given $\Om\in L^{2}(B^{n},so_{m}\otimes\R^{n})$. In the case $n=2$, Rivi\`ere \cite[Theorem I.4]{Riviere-2007} succeeded to prove the existence of $A,B$ under a smallness assumption on $\|\Om\|_{L^2}$. However, as noticed in \cite{Riviere-Struve-2008}, his proof does not work in higher dimensions under the natural Morrey integrability assumption $\Omega\in M^{2,n-2}$ due to the lack of Wente's lemma 
(see \cite{Riviere-Struve-2008,Guo-Xiang-Survey-2022} for more  detailed dicussions).

Motivated by this challenging problem, in this short note, we provide a partial solution in the case $n\ge3$. More precisely, we obtain the following result, which can be viewed as a higher dimensional extension of Rivi\`ere \cite[Theorem 1.4]{Riviere-2007}.

\begin{theorem}\label{thm:existence of A-B} For any $n,m\ge2$, there
exist constants $\ep=\ep(m,n),C=C(m,n)>0$ satisfying the following
property. Suppose $\Om\in L^{n,2}(B^{n},so_{m}\otimes\wedge^{1}\R^{n})$
with
\[
\|\Om\|_{L^{n,2}(B^{n})}\le\ep.
\]
Then there exist $A\in L^{\wq}\cap W^{1,n,1}(B^{n},Gl_{m})$ and $B\in W^{1,n,2}(B^{n},M_{m}\otimes\wedge^{2}\R^{n})$
such that
\begin{equation}
dA-A\Om=-d^{\ast}B.\label{eq: A-B}
\end{equation}
Moreover, we have
\[
\|\dist(A,SO_{m})\|_{L^{\wq}(B^{n})}+\|dA\|_{L^{n,1}(B^{n})}+\|dB\|_{L^{n,2}(B^{n})}\le C\|\Om\|_{L^{n,2}(B^{n})}.
\]
\end{theorem}

In the above theorem, $L^{p,q}$ and $W^{1,p,q}$ denote Lorentz spaces and Lorentz-Sobolev spaces respectively, see e.g.~\cite{Adams-book,ONeil-1963}. We would like to point out that the $L^{n,2}$-integrability condition on $\Omega$ is somehow strong, as $L^{n,2}\subsetneq L^n\subsetneq L^{n,\infty}\subsetneq M^{2,n-2}$ when $n\geq 3$. On the other hand, it seems to be very difficult to refine Rivi\`ere's technique further so that we may relax $L^{n,2}$, for example, to $L^{n,\infty}$.

As that of Rivi\`ere \cite{Riviere-2007}, an application of Theorem \ref{thm:existence of A-B} together with the harmonic mapping equation in the form \eqref{eq: Riviere} leads to

\begin{corollary}\label{coro:conservation law harmonic map}
Let $u\in W^{1,2}(B^{n},N)$ be a weakly harmonic mapping. If $u$ is smooth in $B_{r}(x_{0})\subset B^{n}$ for some
$r>0$ such that $\|\na u\|_{L^{2}(B_{r}(x_{0}))}\ll\ep$, then there
exist $A\in C^{\wq}(B_{r}(x_{0}),Gl_{m}(\R))$ and $B\in C^{\wq}(B_{r}(x_{0}),M_{m}(\R)\otimes\wedge^{2}\R^{n})$
such that $u$ satisfies the conservation law \eqref{eq: CL of Rivere}
of Rivi\`ere in $B_{r}(x_{0})$.
\end{corollary}

We do not know how to extend Corollary \ref{coro:conservation law harmonic map} to neighborhoods around singular points of weakly harmonic mappings in supercritical dimensions and refer to \cite{Guo-Xiang-Survey-2022} for more discussions about the conservation law.


\section{Proof of Theorem \ref{thm:existence of A-B}}

The proof of Theorem \ref{thm:existence of A-B} follows closely that of Rivi\`ere \cite[Theorem 1.4]{Riviere-2007} with minor modifications.

\textbf{Step 1.} There exist $P\in W^{1,n,2}(B^{n},SO_{m})$
and $\xi\in W^{1,n,2}(B^{n},so_{m}\otimes\wedge^{2}\R^{n})$ such
that
\[
P^{-1}dP+P^{-1}\Om P=d^{\ast}\xi\qquad\text{in }B^{n},
\]
and
\begin{equation}\label{eq: epsilon-gauge}
\|dP\|_{L^{n,2}(B^{n})}+\|d\xi\|_{L^{n,2}(B^{n})}\le C\|\Om\|_{L^{n,2}(B^{n})}\le C\ep.
\end{equation}

The proof for this is completely similar to that of \cite[Lemma A.3]{Riviere-2007} or \cite[Lemma 3.1]{Riviere-Struve-2008}, whereas the idea dates back to Uhlenbeck \cite{Uhlenbeck-1982}. Thus we omit it here.

\textbf{Step 2. }To continue, we first extend $\Om$, $P,\xi$ to
the whole space $\R^{n}$ with compact support in a bounded way. By introducing $\tilde{A}=AP$,
to find $(A,B)$ in problem \eqref{eq: A-B} is equivalent to find
$(\tilde{A},B)$ in the problem
\begin{equation}
d\tilde{A}=\tilde{A}d^{\ast}\xi+d^{\ast}BP.\label{eq: equivalent A-B}
\end{equation}
Here and hereafter, all the equations are considered in the whole
space $\R^{n}$. This leads us to solve the second order problem
\begin{equation}\label{eq: A-B-0}
\begin{cases}
-\De\tilde{A}=d^{\ast}\left(\tilde{A}d^{\ast}\xi+d^{\ast}BP\right)=\ast\left(d\tilde{A}\wedge d\left(\ast\xi\right)+d\left(\ast B\right)\wedge dP\right)\\
dd^{\ast}B=d\left((d\tilde{A}-\tilde{A}d^{\ast}\xi)P^{-1}\right)=d\tilde{A}\wedge dP^{-1}+d(\tilde{A}d^{\ast}\xi P^{-1})
\end{cases}
\end{equation}
where we omitted the precise constant coefficients before each term.
Set $\tilde{A}=\hat{A}+\id$. Then we will try to find $(\hat{A},B)$
by solving the equations
\begin{equation}\label{eq: A-B-1}
\begin{cases}
-\De\hat{A}=\ast\left(d\hat{A}\wedge d\left(\ast\xi\right)+d\left(\ast B\right)\wedge dP\right),\\
-\De B=d\hat{A}\wedge dP^{-1}+\ast d^{\ast}\left(\hat{A}d(\ast\xi)P^{-1}+d(\ast\xi)P^{-1}\right),\\
dB=0.
\end{cases}
\end{equation}

\textbf{Step 3.} Solve problem \eqref{eq: A-B-1} by a fixed point argument.

Consider in $\R^{n}$ the problem
\begin{equation}\label{eq: A-B-2}
\begin{cases}
-\De A=\ast\left(da\wedge d\left(\ast\xi\right)+d\left(\ast b\right)\wedge dP\right),\\
-\De B=da\wedge dP^{-1}+\ast d^{\ast}\left(ad(\ast\xi)P^{-1}+d(\ast\xi)P^{-1}\right),\\
dB=0.
\end{cases}
\end{equation}
Introduce the Banach space
\[
\mathbb{X=}\left\{ (a,b):a\in L^{\wq}\cap W^{1,n,2}(\R^{n},Gl_{m}),b\in W^{1,n,2}(\R^{n},M_{m}\otimes\wedge^{2}\R^{n})\right\}
\]
with the norm
\[
\|(a,b)\|_{\X}=\|a\|_{\wq}+\|da\|_{L^{n,2}}+\|db\|_{L^{n,2}}.
\]

For any $(a,b)\in\X$, there exists a unique $(A,B)\in\X$ satisfying
\eqref{eq: A-B-2}. Set
\[
T(a,b)=(A,B).
\]
We prove in below that $T$ is a contraction mapping on $\X_{1}=\{x\in\X:\|x\|_{\X}\le1\}$.

First note that by H\"older's inequality in Lorentz spaces, we have
\[
da\wedge d\ast\xi,da\wedge dP^{-1}\in L^{\frac{n}{2},1}.
\]
Hence from the first equation of \eqref{eq: A-B-2} we derive $A\in W^{2,n/2,1}(\R^{n})\subset C(\R^{n})$
with estimate
\[
\|A\|_{\wq}\le C_{n}\|dA\|_{L^{n,1}}\le C\|d\xi\|_{L^{n,2}}\|da\|_{L^{n,2}}+C\|dP\|_{L^{n,2}}\|db\|_{L^{n,2}}.
\]
From the last two equations of \eqref{eq: A-B-2}, we derive
\[
\|\na B\|_{L^{n,2}}\le\|dP^{-1}\|_{L^{n,2}}\|a\|_{\wq}+C\left(1+\|a\|_{\wq}\right)\|d\xi\|_{L^{n,2}}.
\]
Thus \eqref{eq: epsilon-gauge} implies that
\[
\|(A,B)\|_{\X}\le C\|\Om\|_{L^{n,2}}(\|(a,b)\|_{\X}+1).
\]
Similarly, for any $(a_{i},b_{i})\in\X$, $i=1,2$, we have
\[
\|T(a_{1},b_{1})-T(a_{2},b_{2})\|_{\X}\le C\ep\|(a_{1},b_{1})-(a_{2},b_{2})\|_{\X}.
\]
Thus, by choosing $\ep\ll1$ such that $C\|\Om\|_{L^{n,2}(\R^{n})}\le C\ep\leq 1/2$,
we find that $T$ is a contraction mapping on $\X_{1}=\{x\in\X:\|x\|_{\X}\le1\}$.
Consequently, by the standard fixed point theorem, there exists a
unique $(\hat{A},B)\in\X_{1}$ such that
\[
T(\hat{A},B)=(\hat{A},B).
\]
That is, $(\hat{A},B)$ solves problem \eqref{eq: A-B-1}. Moreover,
\[
\|(\hat{A},B)\|_{\X}\le C\ep(\|(\hat{A},B)\|_{\X}+1)\le2C\ep.
\]

\textbf{Step 4.} Now let $\tilde{A}=\hat{A}+\id$. Then $(\tilde{A},B)$
solves problem \eqref{eq: A-B-0} with
\[
\|\tilde{A}-\id\|_{\wq}+\|d\tilde{A}\|_{L^{n,2}}+\|dB\|_{L^{n,2}}\le C\|\Om\|_{L^{n,2}}.
\]
We need to verify that $(\tilde{A},B)$ solves \eqref{eq: equivalent A-B}.
By the first equation of \eqref{eq: A-B-0} and the nonlinear Hodge theory, there exists $C\in W^{1,n,2}(B^{n},M_{m}\otimes\wedge^{2}\R^{n})$ with $dC=0$ such that
\[
d\tilde{A}-\left(\tilde{A}d^{\ast}\xi+d^{\ast}BP\right)=d^{\ast}C
\]
Then using the second equation of \eqref{eq: A-B-0} we find that
$d(d^{\ast}CP^{-1})=0.$ This gives $d^{\ast}CP^{-1}=dD$ for some
$D\in W^{1,n,2}$, with
\[
\|dD\|_{L^{n,2}}\le C\|\na C\|_{L^{n,2}}\le C\|\na C\|_{L^{n,1}}.
\]
Hence
\[
\De C=dD\wedge dP^{-1},
\]
from which it follows that
\[
\|\na C\|_{L^{n,1}}\le C\ep\|dD\|_{L^{n,2}}.
\]
Thus, combining the above two estimates together gives
\[
\|\na C\|_{L^{n,1}}\le C\ep\|\na C\|_{L^{n,1}}.
\]
By choosing $\ep$ sufficiently small, we obtain $C\equiv0$. The
proof is complete.

\end{document}